\documentclass{article}
\usepackage[utf8]{inputenc}
\usepackage{amsmath,amsfonts,amssymb,amsthm}
\usepackage{dsfont}
\usepackage{todonotes}
\usepackage{color}
\usepackage{algorithm}
\usepackage{comment}
\usepackage{algpseudocode}
\usepackage{algorithmicx}
\usepackage{bm}
\usepackage[makeroom]{cancel}
\usepackage{standalone}
\usepackage{subcaption}

\usetikzlibrary{bayesnet}
\usetikzlibrary{arrows,positioning} 

\DeclareMathOperator*{\argmax}{arg\,max}

\definecolor{sarandonga}{rgb}{0.65,0.35,0}

\definecolor{urdin}{rgb}{0.4,0,0.6}

\newcommand{\M}[0]{{\mathcal M}}
\newcommand{\R}[0]{{\mathbb R}}
\newcommand{\x}[0]{{\mathbf x}}

\title{Parametrization invariant interpretation of priors and posteriors}
\author{Jesus Cerquides\\
IIIA-CSIC, Campus UAB, Cerdanyola\\ 08193 Barcelona, Spain}

\date{}

\begin{document}

\maketitle

%%%%%%%%%%%%%%%%%%%%%%%%%%%%%% Textclass specific LaTeX commands.
\theoremstyle{plain}
\newtheorem{thm}{\protect\theoremname}
\theoremstyle{definition}
\newtheorem{defn}[thm]{\protect\definitionname}
\theoremstyle{plain}
\newtheorem{prop}[thm]{\protect\propositionname}
\theoremstyle{plain}
\newtheorem{lem}[thm]{\protect\lemmaname}
\theoremstyle{definition}
\newtheorem{example}[thm]{\protect\examplename}

\makeatother
\DeclareRobustCommand{\textgreek}[1]{\leavevmode{\greektext #1}}
\providecommand{\definitionname}{Definition}
\providecommand{\examplename}{Example}
\providecommand{\lemmaname}{Lemma}
\providecommand{\propositionname}{Proposition}
\providecommand{\theoremname}{Theorem}

\begin{abstract}
In this paper we leverage on probability over Riemannian manifolds to rethink the interpretation of priors and posteriors in Bayesian inference. The main mindshift is to move away from the idea that ``a prior distribution establishes a probability distribution over the parameters of our model'' to the idea that ``a prior distribution establishes a probability distribution over probability distributions''. To do that we assume that our probabilistic model is a Riemannian manifold with the Fisher metric. Under this mindset, any distribution over probability distributions should be ``intrinsic'', that is, invariant to the specific parametrization which is selected for the manifold. We exemplify our ideas through a simple analysis of distributions over the manifold of Bernoulli distributions. 

One of the major shortcomings of maximum a posteriori estimates is that they depend on the  parametrization. Based on the understanding developed here, we can define the maximum a posteriori estimate which is independent of the parametrization. 
\end{abstract}

\section{Introduction}

Prior and posterior distributions are fundamental to Bayesian inference. As such we are interested in understanding what lies behind a prior or posterior. That is why it is a common practice to plot posterior distributions so that we can perform qualitative analysis such as identifying where the mode of the posterior is, whether it is unimmodal or how spread it is. However, it is well known that priors and posteriors depend on the parametrization we take. In this paper, we follow the well known approach from Jeffreys\footnote{A good description of his philosophy can be found in Section 2.2 in \cite{kass_selection_1996}.} that inference should be invariant to parametrization \cite{jeffreys_theory_1998}. 

%\section{Learning from coin tosses}

Our toy example is the good old statistical example of making inference about a set of coin tosses. Priors for this setting have been widely studied and here we are not interested in the way of selecting a non-informative prior, which has been a Quixotic problem in Bayesian inference. Assume instead that our prior comes given. In particular, we are told by our analyst-in-chief that 
\begin{quotation}
Assume our coin follows a Bernoulli distribution with parameter $\theta$. Our prior distribution over $\theta$ is a $Beta(\frac{1}{2},\frac{1}{2}).$
\end{quotation} 

Enlightened as we are by our analyst-in-chief, we wonder about the kind of information encoded in that prior and to get some intuition we plot its probability density function (pdf), which for the general $Beta(\alpha,\beta)$ we know to be
$$\frac {x^{\alpha -1}(1-x)^{\beta -1}}{\mathrm {B} (\alpha ,\beta )}. $$ 

\begin{figure}
    \centering
    \includegraphics[scale=0.35]{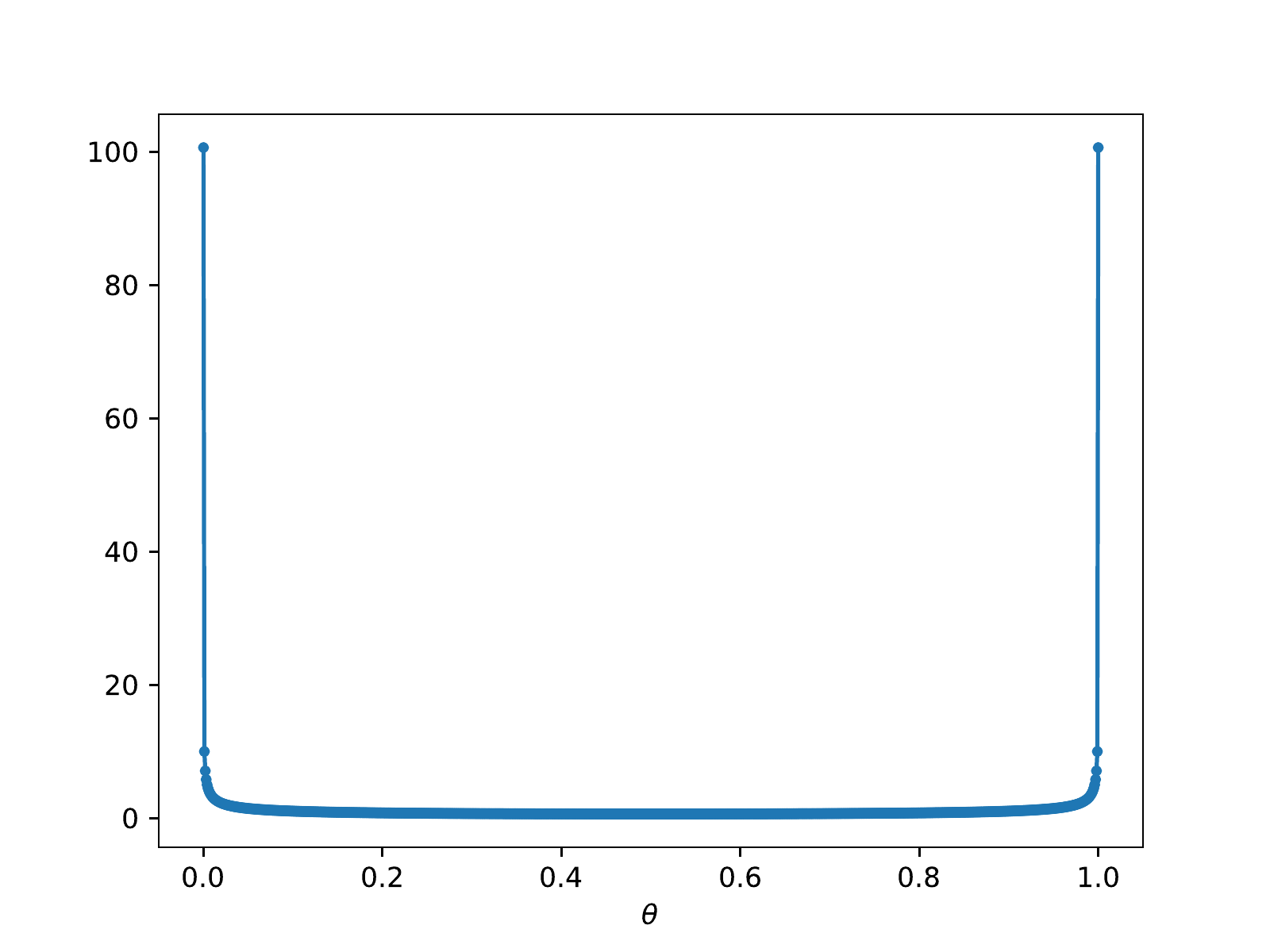}
    \caption{Pdf of a $Beta(\frac{1}{2},\frac{1}{2}).$}
    \label{fig:Beta1}
\end{figure}

After glancing at Figure~\ref{fig:Beta1} we conclude that we are setting a bimodal prior, which has modes at $\theta=0$ and $\theta=1$ and that has minimum density at $\theta=0.5.$

\begin{figure}
    \centering
    \includegraphics[scale=0.35]{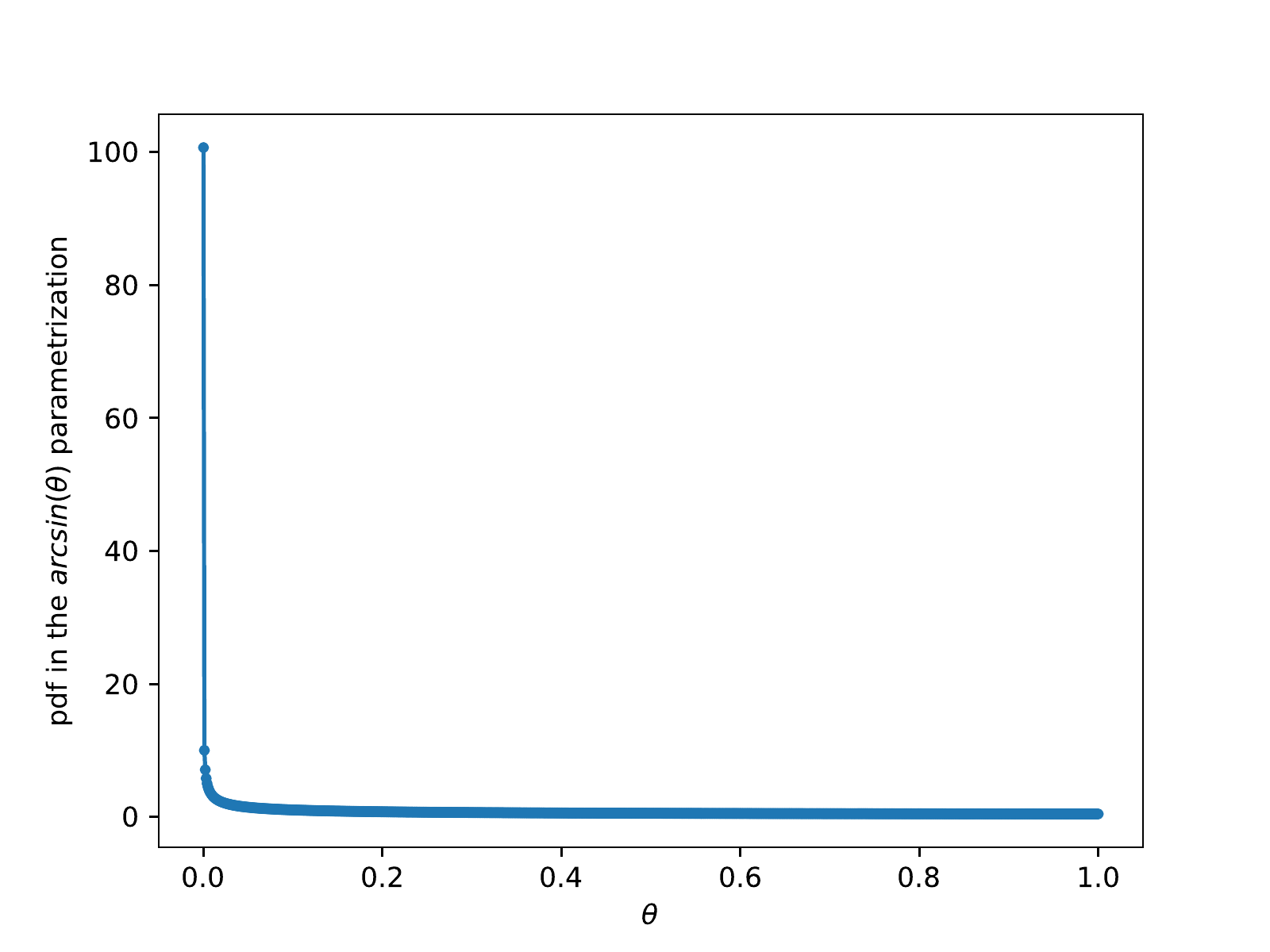}
    \includegraphics[scale=0.35]{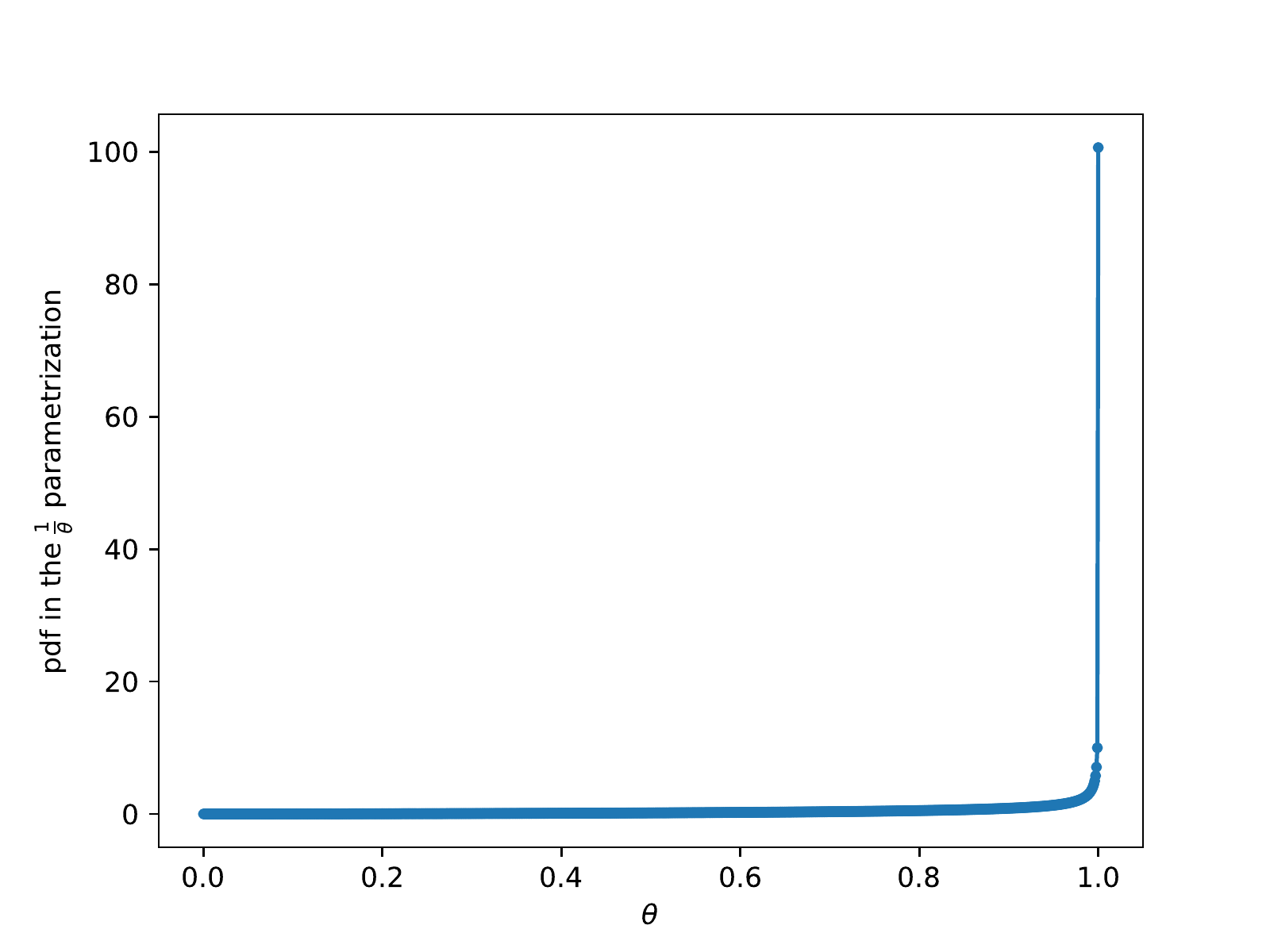}
    \caption{Pdf of a $Beta(\frac{1}{2},\frac{1}{2})$ in the $y=\arcsin(\theta)$ and $y=\frac{1}{\theta}$ parametrizations}
    \label{fig:BetaReparametrizations}
\end{figure}

We go back to our analyst-in-chief and kindly ask him why he thinks the coin is biased. He says, ``you are not looking at it properly, just reparametrize it, try $y=\arcsin \theta$ or $y=\frac{1}{\theta}.$'' So we do, carefully multiplying by the absolute value of the derivative following the rules for change of variable that can be seen in Section~2.6.2 of \cite{murphy_machine_2012}. We are left speechless as we contemplate the density plots with respect to those reparametrizations which can be sen in Figure~\ref{fig:BetaReparametrizations}. 

The same density function, when analyzed in the $y=\arcsin \theta$ parametrization, has a single mode at $\theta=0$, and when analyzed in the $y=\frac{1}{\theta}$ parametrization has a single mode at $\theta=1$. What is happening? Is it possible to really understand what information underlies this density? We will try to provide an answer to this question by using Riemannian geometry. We start by introducing the formal framework required to discuss about probability distributions over probability distributions. 

\section{Probabilities over probabilities\label{sec:ProbsOverProbs}}

Information geometry \cite{amari_information_2016} has shown us that most families of probability distributions can be understood as a Riemannian manifold. In particular in Guy Lebanon words:
\begin{quotation}
A fundamental assumption in the information geometric framework, that manifests itself in many statistical and machine learning applications, is the choice of the Fisher information as the metric that underlies the geometry of probability distributions. The choice
of the Fisher information metric may be motivated in several ways, the strongest of which is Cencov’s characterization theorem (Lemma 11.3 in \cite{cencov_statistical_2000}). ˇ
In his theorem, Cencov proves that the Fisher information metric is the only metric that is invariant under a family of probabilistically meaningful mappings
termed congruent embeddings by a markov morphism. Later on, Campbell \cite{campbell_extended_1986} stripped the proof from its category theory language and extended Cencov’s result to include non-normalized models.
\end{quotation}

Thus, we can work with  probabilities over probabilities by defining random variables which take values in a Riemannian manifold, specifically the one that takes Fisher information as its metric.

Next, we introduce some fundamental definitions required to work with probabilities on a Riemannian manifold. For a more detailed overview of measures and probability see \cite{dudley_real_2002}, of Riemannian manifolds see \cite{jost_riemannian_2011}. Finally, Pennec provides a good overview of probability  on Riemannian manifolds in \cite{pennec_probabilities_2004}.

We start by noting that each manifold $\M$, has an associated $\sigma$-algebra, $\mathcal L_\M$, the Lebesgue $\sigma$-algebra of $\M$ (see section 1, chapter XII in \cite{amann_analysis_2009}). Furthermore, the existence of a metric $g$ induces a measure $\mu_g$ (see section 1.3 in \cite{pennec_probabilities_2004}). The volume of $\M$ is defined as $Vol(\M)=\int_\M 1 d\mu_g.$

\begin{defn}
Let $(\Omega,{\mathcal F},P)$ be a probability space and $(\M, g)$ be a Riemannian manifold. A random variable\footnote{Referred to as a random primitive in \cite{pennec_probabilities_2004}.}  $\mathbf{x}$ taking values in $\M$ is a measurable function from $\Omega$ to $\M.$
Furthermore, we say that $\mathbf{x}$ has a probability density function (p.d.f.) $p_{\mathbf{x}}$ (real, positive, and integrable function) if:
\begin{center}
$\forall {\mathcal X}\in {\mathcal L_\M} \hspace{0.3cm} P(\mathbf{x} \in {\mathcal X} )=\int_{\mathcal X}p_{\mathbf x}d\mu_g$, \hspace{0.5cm} and \hspace{0.5cm} $P(\mathcal M)=1.$
\end{center}
\end{defn}
We would like to highlight that the density function $p_\x$ is intrinsic to the manifold. If $x'=\pi(x)$ is a chart of the manifold defined almost everywhere, we obtain a random vector $\x'=\pi(\x)$. The expression of $p_\x$ in this parametrization is 
\[
p_{\x'}(x')=p_\x(\pi^{-1}(x')).
\]
Let $f:\M\rightarrow \R$ be a real function on $\M$. We define the expectation of $f$ under $\x$ as
\[
\mathbb{E}[f(\x)]=\int_{x}f(x)p_{\x}(x)d{\mu_g}
\]
We have to be careful when computing $\mathbb{E}[f(\x)]$ so that we do it independently of the parametrization. We have to use the fact that $\int_{x}f(x)p_{\x}(x)d{\mu_g}=\int_{x'}f(\pi^{-1}(x'))p_{\x'}(x')\sqrt{\mid G(x')\mid} dx', $
where $G(x')$ is the Fisher matrix at $x'$ in the parametrization $\pi$. Hence, 
\[
\mathbb{E}[f(\x)] =  \int_{x'}f(\pi^{-1}(x'))\rho_{\x'}(x') dx'.
\]
where 
\begin{equation}
\rho_{\x'}(x') = p_{\x'}(x')\sqrt{\mid G(x')\mid} = p_\x(\pi^{-1}(x')) \sqrt{\mid G(x')\mid}
\label{eq:PdfForIntegration}
\end{equation}
is the {\bf expression of $p_\x$ in the parametrization for integration purposes}, that is, its expression with respect to the Lebesgue measure $dx'$ instead of $d\mu_g.$ 

We note that $\rho_{\x'}$ depends on the chart used whereas $p_\x$ is intrinsic to the manifold.

\section{The Riemannian manifold of Bernoulli distributions}

\begin{figure}
    \centering
    \includegraphics[scale=0.35]{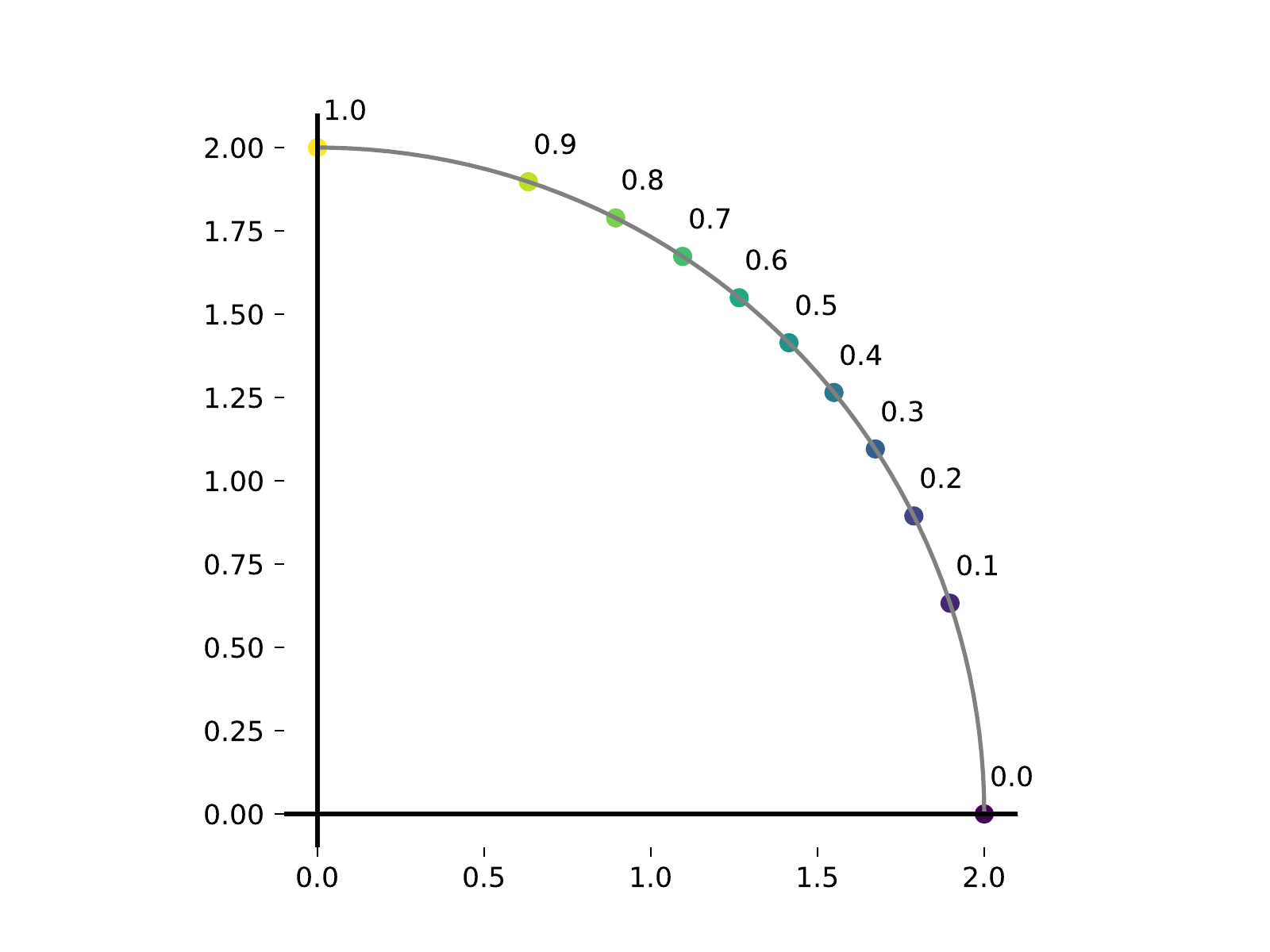}
    \caption{The manifold of Bernoulli distributions with its Fisher metric embedded in ${\mathbb R}^2.$}
    \label{fig:BernoulliManifold}
\end{figure}

In this section we will try to show how the Riemannian geometry of the space of distributions can help us get an understanding of priors and posteriors. 
Our main idea is that both priors and posteriors are probabilities over probability distributions. As such, under reasonable assumptions, they can be understood as probabilities over a Riemannian manifold. 
It is well known that there is an isometry between the set of Bernoulli distributions and the positive quadrant  of the circunference of radius 2 (see section 7.4.2. in \cite{kass_geometrical_1997}). Figure~\ref{fig:BernoulliManifold} shows the embedding of this Riemannian manifold in ${\mathbb R}^2.$ We can observe by the eye that the distance between $0.5$ and $0.6$ is smaller than that between $0.0$ and $0.1.$

\begin{figure}
    \centering
    \includegraphics[scale=0.5]{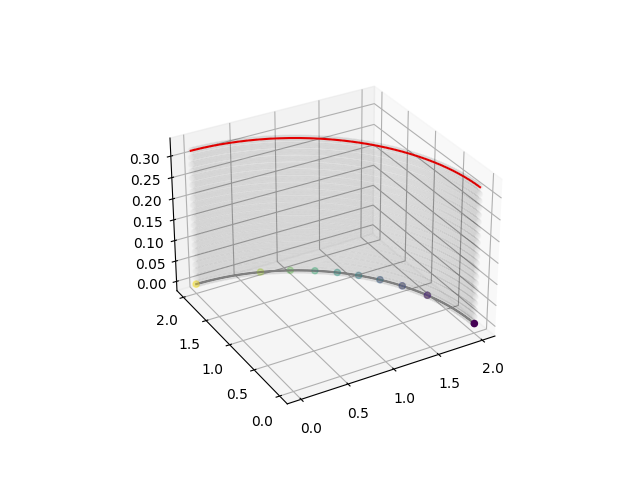}
    \caption{Visualization of the Beta prior as a function over the manifold of Bernoulli distributions}
    \label{fig:Beta3D}
\end{figure}

Now that we know our manifold looks like, we can go back to the problem of understanding the meaning of our prior. Recall that we knew that the prior could be expressed in the $\theta$ parametrization as a $Beta(\frac{1}{2},\frac{1}{2}).$ We can analyze this prior as a random variable (say $\x$) over the manifold of Bernoulli distributions. Note that we use the pdf of the $Beta(\frac{1}{2},\frac{1}{2})$ to compute probabilities by means of integration. Thus, that pdf is just $\rho_{\x'}$, the expression for integration purposes in the $\theta$ parameterization of the probability distribution  $p_\x$ (as presented in Equation~\ref{eq:PdfForIntegration}). Thus, in order to properly understand it we need to find the density function $p_\x$, which is intrinsic to the manifold and hence invariant to reparametrization. From Equation~\ref{eq:PdfForIntegration} it is easy to see that 
\begin{equation}
    p_{\x'}(x')=\frac{\rho_{\x'}(x')}{\sqrt{\mid G(x')\mid }}.\label{eq:ProposedTransformation}
\end{equation}
In our case, 
\begin{equation}
    \rho_{\x'}(\theta)=
    \frac{\theta^{-1/2}(1-\theta)^{-1/2}}
        {B(\frac{1}{2},\frac{1}{2})}
    =\frac{\theta^{-1/2}(1-\theta)^{-1/2}}{\pi},
\end{equation} and the Fisher metric is $G(\theta)=\frac{1}{\theta(1-\theta)}.$ Hence, 
$p_\x(\theta)=\frac{1}{\pi}$, that is, our analyst-in-chief was assuming a uniform distribution over our manifold. To get a geometrical understanding, we can interpret $p_\x(\theta)$ as the height of each of the points in the manifold.  We can visualize this in Figure~\ref{fig:Beta3D}. The height of the red line shows the value of the function $p_\x(\theta)$ as a function of the manifold. The  grey surface that connects the red line with the manifold has area $1$, since the length of the arc is $\pi$, and its height is constant (equal to $\frac{1}{\pi}$). If we want to compute the probability that this distribution assigns to an event, such as for example the probability that $0\leq\theta\leq0.1$, we only need to compute the length of the interval and multiply it by the height of the red line (that is, divide it by $\pi$). 

The same idea that works for Bernoulli distributions works for the many probabilistic models that can be understood as Riemannian manifolds, which includes every member of the exponential family. The main idea to always keep in mind is that the Fisher metric defines not only a distance (the Fisher-Rao distance), but also a volume (a measure, in the measure theoretic sense) on the set of probability distributions. This measure can be used to define the integral of a real function on the manifold. Thus, a probability distribution over the manifold has a unique representation as a function from the manifold to the real numbers which integrates to one. We claim that whenever we want to understand where the mode of a pdf is, how spread a pdf is or simply plot the pdf so that we can grasp an intuitive idea of its behaviour, we should never plot $\rho_{\x'}$ but instead $p_{\x'}$, applying the transformation in Equation~\ref{eq:ProposedTransformation}.

\begin{figure}
    \centering
    \includegraphics[scale=0.5]{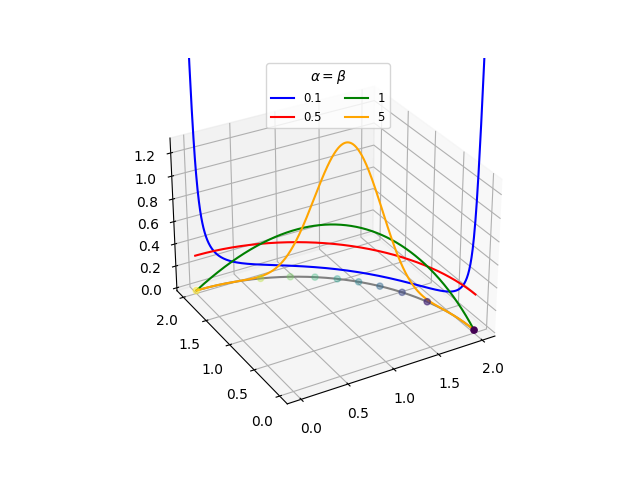}
    \caption{The $p_{\x'}$ density function for different symmetric $Beta$ distributions.}
    \label{fig:SeveralBetas}
\end{figure}

We can use this idea to visualize the knowledge encoded in other Beta distributions. In general, if we assume a $Beta(\alpha,\beta)$ as prior (or we obtain such a pdf as posterior), and want to visualize it or compute its mode, we have to use
\begin{align*} 
p_{\x'}(\theta;\alpha,\beta) &=\frac{\rho_{\x'}(\theta;\alpha,\beta)}{\sqrt{\mid G(\theta) \mid}}\\
&=\frac{\theta^{\alpha-1}(1-\theta)^{\beta-1}}
        {B(\alpha,\beta)}\sqrt{\theta}\sqrt{1-\theta}  \\
&= \frac{\theta^{\alpha-0.5}(1-\theta)^{\beta-0.5}} {B(\alpha,\beta)}.
\end{align*}
Figure~\ref{fig:SeveralBetas} depicts $p_{\x'}$ for a symmetric $Beta$ distribution (that is, $\alpha=\beta$) and for different values of $\alpha$. Specifically, we see that when $\alpha>0.5$ it is a unimodal distribution with mode at $0.5$, and for $\alpha<0.5$ it is a bimodal distribution with modes at $0$ and $1$.  This provides a different understanding from that we get by plotting  $\rho_{\x'}$. Instead, $\rho_{\x'}$, the usual pdf of the $Beta$ function without the correction, shows the change of behaviour from unimodal to multimodal at $\alpha=1$, instead of at $\alpha=0.5.$

\begin{figure}
    \centering
    \includegraphics[scale=0.5]{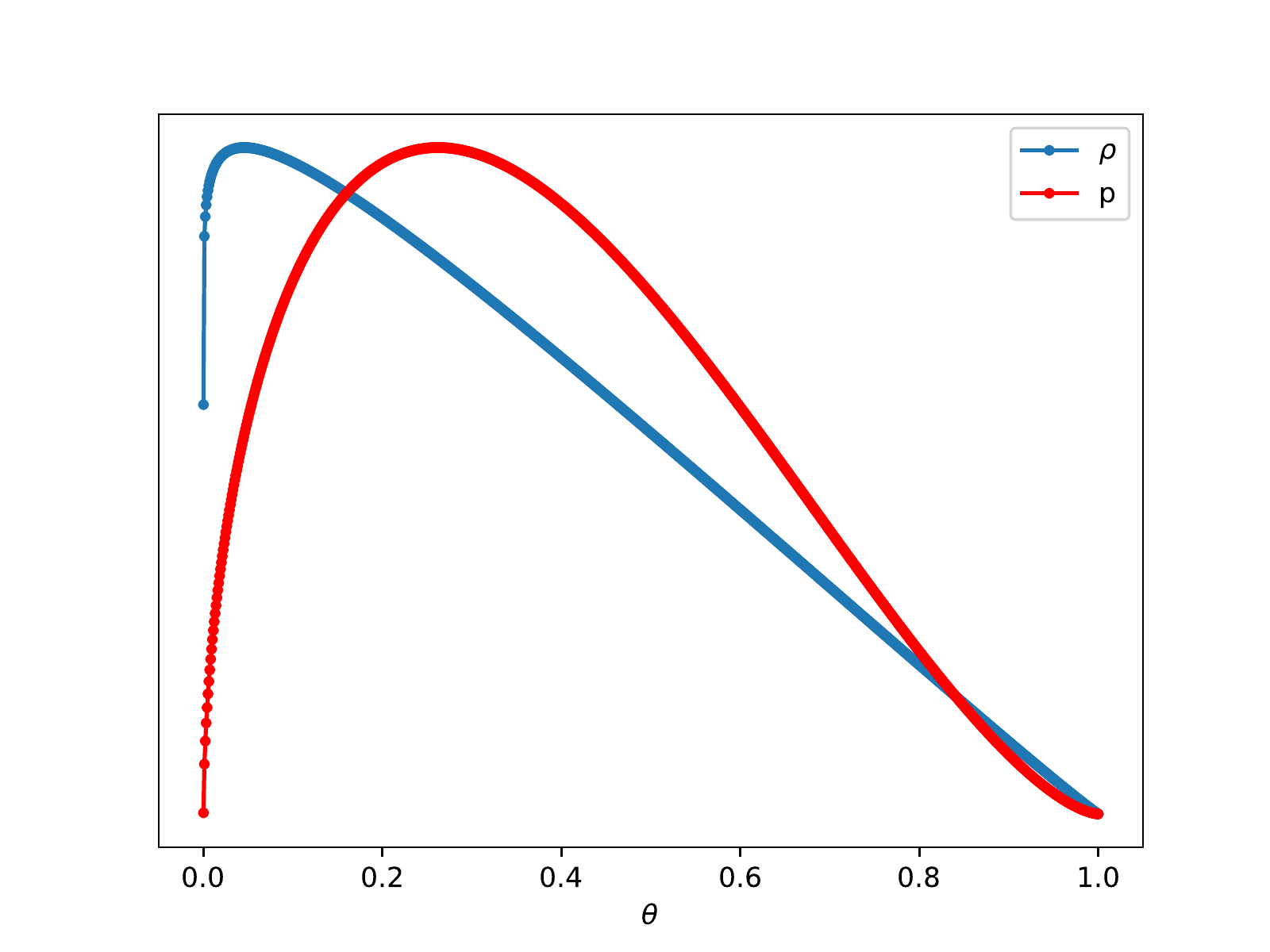}
    \caption{Comparison of $\rho_{\x'}$ and $p_{\x'}$ for $Beta(1.05,2.05).$}
    \label{fig:rho_vs_p}
\end{figure}

Also, when the $Beta$ distribution is asymmetric, the value of the mode of $p_{\x'}$ and that of $\rho_{\x'}$ is different. In Figure~\ref{fig:rho_vs_p}, we can see that the position of the mode is shifted very significantly from around $0.045$ if we compute it according to $\rho_{\x'}$ to $0.26$ if we compute it according to $p_{\x'}.$

\section{Parametrization invariant maximum a posteriori}

The maximum a posteriori (MAP) estimate is a commonly used Bayesian point estimator, defined as the mode\footnote{For simplicity the discussion assumes a single mode. The extension to a set of modes is direct.} of the posterior (see section 4.1.2 in \cite{robert_bayesian_2007}). 

Our coin example has shown us that the computation of the mode depends on the parametrization. In our terminology the usually defined MAP estimate can be expressed as  $$x'_{MAP}=\argmax_{x'}\rho_{\x'}(x').$$
The fact that the MAP estimate is parametrization dependent has brought it criticisms by the community, arguing that it is not a proper Bayesian measure (see for example \cite{bob_carpenter_impact_2016}). 
However, the same Riemannian geometry concepts that allowed us to identify a parametrization invariant notion of pdf ($p_{\x}$) allow us to define a parametrization invariant MAP estimate:
$$x_{MAPI}=\argmax_{x}p_{\x}(x).$$ 
Now, if we need to rely on coordinates for the calculation we can use $p_{\x'}$ and compute 
$$x'_{MAPI}=\argmax_{x'}p_{\x'}(x')$$ in the parametrization of our choice. The probability distribution identified by this procedure will be independent of the parametrization.

\section{Conclusions}
We have shown that by relying on the Riemannian manifold nature of the spaces of probability distributions we can obtain a proper understanding of distributions over probability distributions, by accepting the measure induced by the Fisher metric as a natural measure. The existence of that natural measure can also be used to define the maximum a posteriori estimate independently of the parametrization. 

\section{Acknowledgements} 
Thanks to Borja Sánchez López for the valuable discussions prior to, and along, the writing of this manuscript and for helping me develop a practical understanding of Riemannian geometry. Thanks to Jerónimo Hernández-González for discussions on preliminary versions. This work was partially supported by the projects Crowd4SDG and Humane-AI-net, which have received funding from the European Union’s Horizon 2020 research and innovation program under grant agreements No 872944 and No 952026, respectively. This work was also partially supported by the project CI-SUSTAIN funded by the Spanish Ministry of Science and Innovation (PID2019-104156GB-I00). 

%%%%%%%%%%% The bibliography starts:
\typeout{}
\bibliography{references}

\begin{thebibliography}{10}

\bibitem{amann_analysis_2009}
Herbert Amann and Joachim Escher.
\newblock {\em Analysis {III}}.
\newblock Birkhäuser Basel, 2009.

\bibitem{amari_information_2016}
Sun-ichi Amari.
\newblock {\em Information geometry and its applications}, volume 194.
\newblock Springer, 2016.

\bibitem{bob_carpenter_impact_2016}
{Bob Carpenter}.
\newblock The {Impact} of {Reparameterization} on {Point} {Estimates}.
\newblock In {\em Stan {Case} {Studies}}. April 2016.

\bibitem{campbell_extended_1986}
L.~L. Campbell.
\newblock An extended Čencov characterization of the information metric.
\newblock {\em Proceedings of the American Mathematical Society},
  98(1):135--141, 1986.

\bibitem{cencov_statistical_2000}
N.~N. Cencov.
\newblock {\em Statistical {Decision} {Rules} and {Optimal} {Inference}}.
\newblock American Mathematical Soc., April 2000.
\newblock Google-Books-ID: 63CPCwAAQBAJ.

\bibitem{dudley_real_2002}
R.~M. Dudley.
\newblock {\em Real {Analysis} and {Probability}}.
\newblock Cambridge {Studies} in {Advanced} {Mathematics}. Cambridge University
  Press, Cambridge, 2 edition, 2002.

\bibitem{jeffreys_theory_1998}
Harold Jeffreys.
\newblock {\em The {Theory} of {Probability}}.
\newblock OUP Oxford, August 1998.

\bibitem{jost_riemannian_2011}
Jürgen Jost.
\newblock {\em Riemannian {Geometry} and {Geometric} {Analysis}}.
\newblock Springer, 2011.

\bibitem{kass_geometrical_1997}
Robert~E. Kass and Paul~W. Vos.
\newblock {\em Geometrical foundations of assimptotic inference}.
\newblock Wiley-Interscience, 1997.

\bibitem{kass_selection_1996}
Robert~E. Kass and Larry Wasserman.
\newblock The {Selection} of {Prior} {Distributions} by {Formal} {Rules}.
\newblock {\em Journal of the American Statistical Association},
  91(435):1343--1370, September 1996.
\newblock Publisher: Taylor \& Francis \_eprint:
  https://www.tandfonline.com/doi/pdf/10.1080/01621459.1996.10477003.

\bibitem{murphy_machine_2012}
Kevin~P. Murphy.
\newblock {\em Machine learning: {A} probabilistic perspective}, volume~28.
\newblock MIT Press, 2012.

\bibitem{pennec_probabilities_2004}
Xavier Pennec.
\newblock Probabilities and {Statistics} on {Riemannian} {Manifolds} : {A}
  {Geometric} approach.
\newblock Technical Report RR-5093, INRIA, January 2004.

\bibitem{robert_bayesian_2007}
Christian Robert.
\newblock {\em The {Bayesian} {Choice}: {From} {Decision}-{Theoretic}
  {Foundations} to {Computational} {Implementation}}.
\newblock Springer {Texts} in {Statistics}. Springer-Verlag, New York, 2
  edition, 2007.

\end{thebibliography}
%\hbadness=5000
\bibliographystyle{plain}

\end{document}